\documentclass{amsproc}
\usepackage{amssymb}
\usepackage{amsmath}
\usepackage{amsfonts}
\theoremstyle{plain}

\newtheorem{corollary}{Corollary}
\newtheorem{definition}{Definition}

\newtheorem{theorem}{Theorem}

\begin{document}
\title[Composition Operators]{A note on composition operators on Hardy spaces of the polydisk}
\author{Turgay Bayraktar}
\address{Department of Mathematics\\
Middle East Technical University\\06531 Ankara, Turkey}
\email{turgay@math.metu.edu.tr}
\date{December, 2005}
\subjclass[2000]{Primary 47B33, 47B38}
\keywords{composition operators, Hardy spaces, weakly compact operators}
\thanks{This research has been supported by The Scientific Technological Research Concil of Turkey T\"{U}B\.{I}TAK}

\begin{abstract}
\  In this work, we prove that weak compactness of composition operator on $H^{1}(U^{n})$\ coincides with its compactness. We also characterize bounded and compact composition operators on $H^{1}(U^{n}).$\  
\end{abstract}

\maketitle

\begin{center}
\section{Introduction}
\end{center}
\begin{flushleft}
\ \ Let $\Omega $ be a bounded domain in $\mathbb{C}^{n}$ and $\phi :\Omega
\rightarrow \Omega $ be a holomorphic map. The linear composition operator  
induced by the symbol $\phi $ is defined by 
\[
C_{\phi }f=f\circ \phi \ , \ \ f\in H(\Omega).
\]
\ \ Boundedness and compactness of composition operators on $H^{p}$ spaces were
studied by various authors \cite{b10},\cite{b9}, \cite{b6} and \cite{b8}. As a consequence of Littlewood's
subordination principle \cite{b4}, every holomorphic self-map $\phi $ of the unit
disc induces a bounded composition operator on $H^{p}(\mathbb{D}).$ However,
boundedness of composition operators is not automatic in several variables
(see \cite{b9} and \cite{b6}). In \cite{b14} J. Shapiro characterized compact composition
operators in terms of Nevanlinna counting function that will not be defined
in this paper. MacCluer \cite{b9} provided a different characterization of bounded
(compact) composition operators on Hardy spaces of the unit ball  
$H^{p}(B^{n})$ in terms of Carleson measures for $0<p<\infty .$ In the
polydisk case, Jafari \cite{b6} obtained the analogues results on $H^{p}(U^{n})$
for $p>1.$ However, the case $p=1$ remained open in his work \cite{b6}.\newline

\ \ Later, Sarason \cite{b13} showed that weak compactness of composition operator on 
$H^{1}(\mathbb{D})$ implies its compactness. More recently, Li and Russo \cite{b8}
generalized Sarason's result to strongly pseudoconvex domains. However, the
polydisk case remained open. Indeed, it was posed as an open problem in \cite{b12}.

\ \ In the first part of this work, we extend Jafari's result on bounded
(compact) composition operators on $H^{p}(U^{n})$ \ to $p=1$ \ by making use
of a weak factorization theorem \cite{b5} for $H^{1}(U^{n}).$ Finally, 
we prove that weak compactness of composition operator
on $H^{1}(U^{n})$ coincides with its compactness. The argument that we
present is rather general and unifying.

\section*{Preliminaries}
We denote the \emph{open unit polydisk} by
\begin{equation*}
U^{n}=\{z\in\mathbb{C}^{n}:\mid z_{i}\mid< 1\text{\ for all i=1,...,n}\} 
\end{equation*}
and the \emph{distinguished boundary} of the polydisk $U^{n}$ by
\begin{equation*}
T^{n}=\{z\in\mathbb{C}^{n}:\mid z_{i}\mid=1\text{ \ for all i=1,...,n }\}
\end{equation*}
The \emph{Hardy spaces} are defined by 
\begin{align*}
H^{p}(U^{n}) & =\{f\in H(U^{n}):\parallel f\parallel_{p}<\infty\}
\text{ \ \ where}\ 0<p<\infty\ and \\
& \parallel f\parallel_{p}\doteq\sup_{0<r<1}((2\pi)^{-n}\int_{0}^{2\pi
}...\int_{0}^{2\pi}\left|
f(re^{i\theta_{1}},...,re^{i\theta_{n}})\right|^{p}d\theta_{1}...d\theta_{n})^{1/p}\text{\ }
\end{align*}
and 
\begin{equation*}
H^{\infty}(U^{n})=\{f\in H(U^{n}):\parallel f\parallel_{\infty}\doteq
\sup_{z\in U^{n}}\left| f(z)\right|<\infty\} 
\end{equation*}
Let $f:U^{n}\rightarrow \mathbb{C}$\ be a holomorphic map. Then the \textit{radial limits} of $f$\ are given by
\begin{equation*}
f^{\ast }(\zeta)=\lim_{r\rightarrow 1^{-}}f(r\zeta),\ \ \zeta \in T^{n}
\end{equation*}
whenever the limit exists. In fact, we have 
\begin{equation*}
\parallel f\parallel_{p}^{p}=\int_{0}^{2\pi}...\int_{0}^{2\pi}\left|f^{*}(\zeta)\right|^{p}dm_{n}(\zeta)
\end{equation*}
where $m_{n}=(2\pi)^{-n}d\theta_{1}...d\theta_{n}$\ denotes the normalized Lebesgue area measure on $T^{n}$ (\textit{see \cite{b11}} for details). That is , $H^{1}(U^{n})$\ is isometrically embedded into $L^{1}(m_{n},T^{n})$. 
\newline
For $f\in L^{p}(T^{n})$, $1\leq p<\infty$\ we denote the Poisson integral of $f$\ by $P[f]$.
\newline  
Let I be an interval on T of \textit{length} $\left|I\right|=\delta$\ centered at $\zeta=e^{i\theta_{0}}$, S(I) is defined by
\begin{equation*}
S(I)=\{z\in U:z=re^{i\theta}, 1-\delta< r <1\ and \ \left|\theta-\theta_{0}\right|<\frac{\delta}{2}\}.
\end{equation*}
If $R=I_{1}\times I_{2}\times...\times I_{n}\subset T^{n}$, with $I_{j}$\ intervals having length $\delta_{j}$ and centers $\zeta_{j}$, then S(R) is given by
\begin{equation*}
S(R)=S(I_{1})\times S(I_{2})\times...\times S(I_{n})
\end{equation*}
and
\begin{equation*}
m_{n}(R)=\left|I_{1}\right|\times...\times\left|I_{n}\right|
\end{equation*}
for every rectangle $R=I_{1}\times..\times I_{n} \subset T^{n}.$\newline
Let V is an open set in $T^{n}$,\ S(V) be defined by 
\begin{equation*}
S(V)=\bigcup\{S(R): R=I_{1}\times...\times I_{n}\subset V\}.
\end{equation*} 

Following \cite {b1}, a positive Borel measure $\mu$\ on\ $U^{n}$\ is said to be a \emph{Carleson measure} if there exists a constant $C>0$\ such that 
\begin{equation}
\mu(S(V))\leq Cm_{n}(V)\ 
\end{equation}   
for all connected open sets\ $V\subset T^{n}$.\newline The measure $\mu$\ is said to be a \emph{compact-Carleson measure} \cite{b7} if 
\begin{equation}
\lim_{\delta\rightarrow 0^{+}}(\sup\{\frac{\mu(S(V))}{m_{n}(V)}: V\subset T^{n}\ open\ connected,\ m_{n}(V)<\delta\})=0.	
\end{equation}
It is a well-known fact that every function $f\in H^{1}(U)$\ can be expressed as a product of two functions, $f=gh$\ where $g\ and\ h\in H^{2}(U).$\ This factorization theorem fails in several variables (\textit{See \cite[Sec. 4.2]{b11}} ).  However, there exists a weak factorization theorem (\textit{see} \cite{b5}). 
\begin{theorem}
\ If \ $f \in H^{1}(U^{n})$\ then there exists\ $g_{j},\ h_{j}\in H^{2}(U^{n}),$\ such that
\begin{equation*}
f=\sum_{j=1}^{\infty}g_{j}h_{j}
\end{equation*}
and $\sum_{j=1}^{\infty} \parallel g_{j}\parallel_{H^{2}} \parallel h_{j}\parallel_{H^{2}} \leq C\parallel f\parallel_{H^{1}}$ for some $C>0$.
\end{theorem}

\section{Results}

\subsection{Boundedness and compactness of $C_{\phi}$}

\bigskip
\begin{flushleft}
In this part, we will observe boundedness and compactness of composition operators on $H^{1}(U^{n})$.\newline
In \cite{b1}, Chang obtained a characterization of Carleson measures for $1<p<\infty$. The following theorem asserts that Chang's result extends to $p=1$.
\begin{theorem}
Suppose that $\mu$\ is a positive Borel measure on $\overline{U^{n}}$\ and $1\leq p<\infty$\ then the following conditions are equivalent
\begin{itemize}
	\item [(i)] There\ exists\ a\ constant\ $C>0$\ such\ that
	\begin{equation*}
	\int_{\overline{U^{n}}}\left|P[f(z)]\right|^{p} d\mu(z)\leq C \int_{T^{n}}\left|f^{*}(\zeta)\right|^{p}dm_{n}(\zeta)
	\end{equation*}
	for all\ $f\in H^{p}(U^{n}).$   
  
  \item [(ii)] There exists a constant $C>0$\ such that 
\begin{equation*}
\mu(S(V))\leq Cm_{n}(V)\ 
\end{equation*}   
for all connected open sets\ $V\subset T^{n}$.
\end{itemize}
  
\end{theorem}

\begin{proof}
It is enough to prove for the case p=1.\newline
Suppose that condition (i) holds. Since $H^{2}(U^{n})\subset H^{1}(U^{n})$\ and $H^{2}(U^{n})$\ is a complemented subspace of $L^{2}(T^{n})$, by closed graph theorem we have
\begin{equation*}
\int_{\overline{U^{n}}}\left|P[f]\right|d\mu \leq C\int_{T^{n}}\left|f\right|dm_{n}\  for\ all\ f\in L^{2}(T^{n}).
\end{equation*} 
Now $(ii)$\ follows by taking $f\equiv \chi_{V}.$

  To show the converse implication, let us assume that $(ii)$\ holds. Then we have $H^{2}(U^{n})\subset L^2(\mu,\overline{U^{n}})$ (\textit{see \cite{b1}}). Let $f\in H^{1}(U^{n})$. Then there exists $\{g_{j}\}$\ and $\{h_{j}\}$\ in $H^{2}(U^{n})$\ such that $f=\sum_{j=1}^{\infty}g_{j}h_{j}$. \newline We consider the decomposition $\mu=\mu_{1}+\mu_{2}$\ where $\mu_{1}$\ is the restriction of $\mu$ to $U^{n}.$\ Since $\mu$\ is a Carleson measure, $supp(\mu_{2}) \subset T^{n}.$\newline
Note that $\sum_{j=1}^{\infty}g_{j}h_{j}$\ converges to $f$\ pointwise in $U^{n}$\ and it is easy to show that radial limits of $\sum_{j=1}^{\infty}g_{j}h_{j}$\ converge a.e. pointwise to radial limits of $f$\ on $T^{n}$. Then, first by Fatou's lemma and then by H\"{o}lder's inequality, we have 
\begin{align*}
\int_{\overline{U^{n}}}\left|f\right|d\mu &\leq \sum_{j}\int_{\overline{U^{n}}}\left|g_{j}\right|\left|h_{j}\right|d\mu \leq \sum_{j}(\int_{\overline{U^{n}}}\left|g_{j}\right|^{2} d\mu )^{\frac{1}{2}}(\int_{\overline{U^{n}}}\left|h_{j}\right|^{2} d\mu )^{\frac{1}{2}} \\
&\leq C_{1} \sum_{j} \parallel g_{j}\parallel_{H^{2}} \parallel h_{j}\parallel_{H^{2}} \\
&\leq C \parallel f \parallel_{H^{1}}
\end{align*}
\end{proof}

\begin{theorem}
Let $\mu$\ be a positive Borel measure on $\overline{U^{n}}$\ and assume that $\mu$\ is a Carleson measure. Let I be the natural embedding sending $H^{p}(U^{n})$\ to $L^{p}(\mu,\overline{U^{n}})$\ , $1\leq p<\infty.$\ Then the following two conditions are equivalent
\begin{itemize}
	\item [(i)] I is compact
	\item [(ii)]$\mu$ is a compact-Carleson measure, i.e.
	\begin{equation*}
	\lim_{\delta\rightarrow 0^{+}} (\sup\{\frac{\mu(S(V))}{m_{n}(V)}: V\subset T^{n}, m_{n}(V)<\delta\})=0
	\end{equation*}
\end{itemize}
\end{theorem}
\begin{proof}
The same proof as in \cite[Theorem 2.3]{b7} applies for $p=1$.
\end{proof}
Let $\mu$\ denote the pull back measure  
\begin{equation*}
\mu (E)=m_{n}((\phi ^{\ast })^{-1}(E))\text{ \ for\ }E\subset \overline{U^{n}} 
\end{equation*}
where $dm_{n}$\ is the normalized Lebesgue area measure on $T^{n}.$\ Then we have the following change of variable formula 
\begin{equation*}
\int\limits_{T^{n}}\mid (f\circ \phi)^{*} \mid^{p} dm_{n}=\int\limits_{T^{n}}\mid f^{*}\circ \phi^{*} \mid^{p} dm_{n}=\int\limits_{\overline{U^{n}}}\mid f\mid^{p} d\mu
\end{equation*}
for every $f\in H^{p}(U^{n})\cap C(\overline{U^{n}})$\ \emph{(see \cite{b2} and \cite{b6} for details).}\newline 
As a consequence of Theorem 2 and Theorem 3 and using a density argument, we have the following corollary:
\begin{corollary}
Let $\phi:U^{n}\rightarrow U^{n}$\ be holomorphic and $1\leq p<\infty$. $C_{\phi}$\ is bounded (resp. compact) on $H^{p}(U^{n})$\  if and only if  $m_{n}((\phi^{*})^{-1})$\ is a Carleson measure (resp. compact Carleson measure).
\end{corollary}
\end{flushleft}
\subsection{Compactness and Weak Compactness}
\bigskip

\begin{flushleft}
In this section we will observe the relationship between compactness and weak compactness of composition operator.\newline
The linear operator on a Banach space E is said to be compact (respectively weakly compact) if it maps unit ball of E to a relatively compact (respectively relatively weakly compact) subset of E. \\
In the sequel we will need the following notion of uniformly integrable subsets of $L^{1}(\mu,X)$.
\begin{definition}
Suppose that $\mu$\ is a positive Borel measure on a topological space X. A bounded set $\Lambda\subset L^{1}(\mu,X)$\ is said to be uniformly integrable if  $\forall \alpha>0,\ \exists \beta>0$\ such that
\begin{equation*}
\int_{E}\left|f\right|d\mu<\alpha \ whenever\ \mu(E)<\beta\ and\ f\in \Lambda.
\end{equation*}      
\end{definition}
A well known result of N. Dunford \cite[p. 294]{b3} asserts that for a finite measure $\mu$, a subset $\Lambda$\ of $ L^{1}(\mu,X)$\ is weakly relatively compact if and only if $\Lambda$\ is uniformly integrable . 

\begin{theorem}
Let $\phi :U^{n}\rightarrow U^{n}$ \ be holomorphic map. Then the following are
equivalent:

\begin{itemize}
\item [(i)] $C_{\phi }:H^{1}(U^{n})\rightarrow H^{1}(U^{n})$ is compact

\item [(ii)] $C_{\phi }:H^{1}(U^{n})\rightarrow H^{1}(U^{n})$ is weakly compact

\item [(iii)] $m_{n}((\phi^{*})^{-1})$\ is a compact-Carleson measure, i.e.
\begin{equation*}
\lim_{\delta\rightarrow 0^{+}}(\sup\{\frac{\mu(S(V))}{m_{n}(V)}: V\subset T^{n}\ open\ connected,\ m_{n}(V)<\delta\})=0	
\end{equation*}
\end{itemize}
\end{theorem}

\begin{proof}
$(i\Rightarrow ii)$ Obvious

$(ii\Rightarrow iii)$ Since $C_{\phi }:H^{1}(U^{n})\rightarrow H^{1}(U^{n})$ is weakly compact, $C_{\phi}$\ is bounded. Then by Theorem 2, $\mu=m_{n}((\phi^{*})^{-1})$\ is a Carleson measure. That is, the identity 
\begin{equation*}
I:H^{1}(U^{n})\rightarrow L^{1}(\mu,\overline{U^{n}}) 
\end{equation*}
\begin{equation*}
I(f)=f
\end{equation*}
is bounded.
Moreover, we have
\begin{equation}
\parallel I(f)\parallel=\parallel f\parallel _{L^{1}}=\parallel C_{\phi }f\parallel_{H^{1}}. 
\end{equation}
Let $B$\ denotes the unit ball of $H^{1}(U^{n})$. Using weak compactness of $C_{\phi}$\ and the isometric embedding of $H^{1}(U^{n})$\ into $L^{1}(m_{n},T^{n})$, one gets $C_{\phi}(B)$\ is uniformly integrable in $L^{1}(m_{n},T^{n})$. Then applying the change of variable formula in (3), we get I(B) is uniformly integrable in $L^{1}(\mu,\overline{U^{n}})$. Thus, I is weakly compact.\newline
Now, suppose that  condition $(iii)$\ does not hold. Then $\exists $ 
$\varepsilon >0$ , $V_{j}\subset T^{n}$\ such that $lim_{j\rightarrow 0}m_{n}(V_{j})=0$\ and  
\begin{equation*}
\mu(S(V_{j}))>\varepsilon m_{n}(V_{j}). 
\end{equation*}
Then there exists half-open rectangles $R_{j}\subset V_{j}$\ such that $\lim_{j\rightarrow\infty} m_{n}(R_{j})=0$\ and  $\mu(S(R_{j})) \geq \frac{\epsilon}{2}m_{n}(R_{j}).$\ Indeed, we can express $V_{j}=\bigcup_{k=1}^{\infty} R_{k}^{j}$\ for some half-open rectangles $R_{k}^{j}\subset T^{n}$\ with $Int R_{k}^{j}\cap Int R_{l}^{j}=\varnothing$\ whenever $k\neq l$\ and suppose that $\mu (S(R_{k}^{j})) <\frac{\epsilon}{2}m_{n}(R_{k}^{j})$\ then 
\begin{equation*}
\mu(S(V_{j}))\leq \sum_{k}\mu(S(R_{k}^{j}))<\frac{\epsilon}{2}\sum_{k}m_{n}(R_{k}^{j})=\frac{\epsilon}{2}m_{n}(V_{j}) 
\end{equation*}
this contradicts $\mu(S(V_{j}))>\epsilon m_{n}(V_{j}).$ \newline
Define $f_{j}$\ on $U^{n}$\ as
\begin{equation*}
f_{j}(z_{1},..,z_{n})=\prod\limits_{i=1}^{n}(1-\overline{\alpha}_{j_{i}}z_{i})^{-4}
\end{equation*}
where $\alpha_{j_{i}}=(1-\delta_{j_{i}})\zeta_{j_{i}}.$\ Then we have
\begin{equation*}
\parallel f_{j}\parallel_{1}\ \sim (\prod_{i=1}^{n} (1-\left|\alpha_{j_{i}}\right|^{2}))^{-3}\ \sim (\prod_{i=1}^{n}\delta_{j_{i}})^{-3}.
\end{equation*}
If $z\in S(R_{j})$\ 
\begin{align*}  
\left|1-\overline{\alpha_{j_{i}}}z_{i}\right|&=\left|1-(1-\delta_{j_{i}})\overline{\zeta_{j_{i}}}z_{i}\right| \\
&\leq\left|\overline{\zeta_{j_{i}}}(\zeta_{j_{i}}-z_{i})\right|+\left|\delta_{j_{i}}\overline{\zeta_{j_{i}}}z_{i}\right| \\
&\leq 2\delta_{j_{i}}. 
\end{align*}
Then
\begin{equation*}
\left|f_{j}(z)\right|\geq (\frac{1}{16})^{n}(\prod_{i=1}^{n}\delta_{j_{i}})^{-4}\  on\ S(R_{j}).
\end{equation*} 
Let  $g_{j}=f_{j}/\parallel f_{j}\parallel_{H^{1}}.$\ It is enough to show that there
exists a subsequence $g_{j_{k}}$ of  $g_{j}$\ for which the set \{$I(g_{j_{k}})\mid k\in \mathbb{N}$\} is not uniformly integrable.
\newline
Since $\mu=m_{n}((\phi^{*})^{-1})$\ is a Carleson measure $\exists\ M>0$\ such that
\begin{equation*}
\mu(S(R))\leq Mm_{n}(R)\ \ for\ every\ rectangle\ R\subset T^{n}   
\end{equation*}
Fix $\beta>0$, then $\exists\ k_{0}$\ such that
\begin{equation*}
\mu(S(R_{j_{k}}))\leq \beta\ whenever\ k\geq k_{0}. 
\end{equation*}   
Then
\begin{equation*}
\left|f_{j_{k}}(z)\right| \geq (\frac{1}{16})^{n} (\prod_{i}\delta_{j_{k_{i}}})^{-4}\ whenever\ z\in S(R_{j_{k}}).
\end{equation*}
Thus, we have
\begin{align*}
\int_{S(R_{j_{k}})}\left|g_{j_{k}}\right|d\mu &\geq\frac{(\frac{1}{16})^{n}(\prod_{i}\delta_{j_{k_{i}}})^{-4}}{\parallel f_{j_{k}}\parallel_{H^{1}}}\mu(S(R_{j_{k}})) \\
&\geq (\frac{1}{16})^{n}\frac{\epsilon}{2}\frac{(\prod_{i}\delta_{j_{k_{i}}})^{-4}}{\parallel f_{j_{k}}\parallel_{H^{1}}}(\prod_{i}\delta_{j_{k_{i}}}) \\
&\geq \frac{\epsilon}{2.16^{n}} 
\end{align*}
Hence, $I(\{g_{j_{k}}\})$\ is not uniformly integrable. \newline 
$(iii\Rightarrow i)$ Follows from Corollary 1.  
\end{proof}

\end{flushleft}

\end{flushleft}


\begin{thebibliography}{14}

\bibitem{b1}S.-Y.A. Chang, Carleson measure on the bi-disc, Ann. of Math. (2) \textbf{109} (1979), 613-620. 

\bibitem{b2}C.C. Cowen \& B. MacCluer, Composition Operators on Spaces of
Analytic Functions, CRC Press, 1995. 

\bibitem{b3}N. Dunford \& J.T. Schwartz, Linear Operators Part I, Wiley-Intersience, New York, 1958.
 
\bibitem{b4}P. Duren, Theory of H$^{p}$ spaces, Academic Press, 1970. 

\bibitem{b5}S. H. Ferguson \& M. T. Lacey, A characterization of product BMO by commutators, Acta Math., \textbf{189}(2002) 2, 143-160.

\bibitem{b6}F. Jafari, On bounded and compact composition operators in
polydiscs, Canadian J. Math. XLII (1990), 869-889. 

\bibitem{b7}F. Jafari, Carleson measures in Hardy and weighted Bergman spaces of polydiscs, Proc. Amer. Math. Soc.\textbf{112}(1991), 771-781.

\bibitem{b8}S.-Y. Li \& B. Russo, On compactness of composition operators in Hardy spaces of several variables, Proc. Amer. Math. Soc. \textbf{123}(1995), 161-171.

\bibitem{b9}B.MacCluer, Compact composition operators on $H^{p}(B^{N})$, Michigan J. Math. \textbf{32}(1985), 237-248. 

\bibitem{b10}B. MacCluer \& J.H. Shapiro, Angular derivatives and compact
composition operators on the Hardy and Bergman spaces, Canadian J. Math. \textbf{38}
(1986), 878-906. 

\bibitem{b11}W. Rudin, Function theory in polydiscs, Benjamin, New York, 1969.


\bibitem{b12}B. Russo, Holomorphic composition operators in several complex variables, Studies on Composition Operators,  Proceedings of the Rocky Mountain Mathematics Consortium 1996, Contemp. Math. \textbf{213} (1998), 191-212. 

\bibitem{b13}D. Sarason, Weak compactness of holomorphic composition operators
on H$^{1},$\ Functional Analysis and Operator Theory (New Delhi,1990), 75-79
Sipringer-Verlag, Berlin, 1992.

\bibitem{b14}J.H. Shapiro, The essential norm of a composition operator,
Annals Math. \textbf{125}(1987), 375-404.

\end{thebibliography}
\end{document}